\documentclass[11pt]{article}
\usepackage{epic,latexsym,amssymb}
\usepackage{color}
\usepackage{tikz}
\usepackage{amsmath}

\newtheorem{theorem}{Theorem}

\newtheorem{fact}{Fact}
\newtheorem{claim}{Claim}

\newtheorem{corollary}{Corollary}

\newtheorem{conjecture}{Conjecture}
\newtheorem{problem}{Problem}

\newtheorem{definition}{Definition}

\newcommand{\QED}{$\Box$}

\newcommand{\smallqed}{{\tiny ($\Box$)}}

\newcommand{\cG}{{\cal G}}

\newcommand{\w}{{\rm w}}

\newcommand{\pdom}{{\rm pdom}}
\newcommand{\gpr}{\gamma_{\rm pr}}

\newcommand{\proof}{\noindent\textbf{Proof. }}

\newcommand{\1}{ \vspace{0.1cm} }

\let\oldenumerate\enumerate
\renewcommand{\enumerate}{
  \oldenumerate
  \setlength{\itemsep}{0pt}
  \setlength{\parskip}{0pt}
  \setlength{\parsep}{0pt}
}

\def\vertex(#1){\put(#1){\circle*{2}}}
\def\vertexo(#1){\put(#1){\circle{2}}}
\def\vert(#1){\put(#1){\circle*{1.5}}}
\def\verto(#1){\put(#1){\circle{1.5}}}
\def\lab(#1)#2{\put(#1){\makebox(0,0)[c]{#2}}}

\begin{document}

\title{Paired domination in graphs with minimum degree four}

\author{$^{1,2}$Csilla Bujt\'{a}s \, and \, $^3$Michael A. Henning \\ \\
$^1$Faculty of Mathematics and Physics \\
University of Ljubljana\\
Ljubljana, Slovenia \\
\\
$^2$Institute of Mathematics, Physics, and Mechanics \\
Ljubljana, Slovenia \\
\small \tt Email: csilla.bujtas@fmf.uni-lj.si \\
\\
$^3$Department of Mathematics and Applied Mathematics \\
University of Johannesburg \\
Auckland Park, 2006 South Africa\\
\small \tt Email: mahenning@uj.ac.za  \\
}
 	
\date{}
\maketitle

\begin{abstract}
A set $S$ of vertices in a graph $G$ is a paired dominating set if every vertex of $G$ is adjacent to a vertex in $S$ and the subgraph induced by $S$ admits a perfect matching. The minimum cardinality of a paired dominating set of $G$ is the paired domination number $\gpr(G)$ of $G$. We show that if $G$ is a graph of order~$n$ and $\delta(G) \ge 4$, then $\gpr(G) \le \frac{10}{17}n < 0.5883 n$.
\end{abstract}

{\small \textbf{Keywords:} Paired domination; Bounds; Minimum degree four } \\
\indent {\small \textbf{AMS subject classification:} 05C69}

\date{}
\maketitle

\section{Introduction}

A set $S$ of vertices in a graph $G$ is a \emph{dominating set} if every vertex in $V(G) \setminus S$ is adjacent to a vertex in~$S$. The \emph{domination number} $\gamma(G)$ of $G$ is the minimum cardinality of a dominating set of $G$. The study of domination in graphs is a vibrant and growing area of research in graph theory. A thorough treatment of this topic can be found in the recent so-called ``domination books''~\cite{HaHeHe-20,HaHeHe-21,HaHeHe-23,HeYe-book}.
A \emph{paired dominating set}, abbreviated PD-set, of $G$ is a dominating set $S$ of $G$ such that the induced subgraph $G[S]$ contains a perfect matching $M$ (not necessarily induced). With respect to the matching $M$, two vertices joined by an edge of $M$ are \emph{paired} and are called \emph{partners} in $S$. The \emph{paired domination number}, $\gpr(G)$, of $G$ is the minimum cardinality of a PD-set of $G$. A $\gpr$-\emph{set} of $G$ is a PD-set of $G$ of minimum cardinality. Necessarily, the paired domination number of a graph is an even integer. Paired domination in graphs is well studied in the literature, and was first studied by Haynes and Slater~\cite{HaSl-95,HaSl-98} in 1995. Surveys on paired domination in graphs can be found in~\cite{DeHe-14,DeHaHe-20}. In this paper, we study bounds on the paired domination number of a graph with minimum degree at least~$4$.

\subsection{Graph theory notation and terminology}
 	
For notation and graph theory terminology, we generally follow~\cite{HaHeHe-23}. Specifically, let $G$ be a graph with vertex set $V(G)$ and edge set $E(G)$, and of order~$n(G) = |V(G)|$ and size $m(G) = |E(G)|$. If $S$ is a set of vertices in $G$, the graph $G - S$ denotes the graph obtained from $G$ by removing the vertices (and their incident edges) from $S$. If $S = \{v\}$, then we simply write $G - v$ rather than $G - \{v\}$. The subgraph induced by the set $S$ is given by $G[S]$. The \emph{open} (resp., \emph{closed}) \emph{neighborhood}  of the set $S \subseteq V(G)$ is the union of the open (resp., closed) neighborhoods of vertices in $S$, denoted by $N_G(S)$ (resp., $N_G[S]$).

We denote the \emph{degree} of $v$ in $G$ by $\deg_G(v) = |N_G(v)|$. The minimum and maximum degree in $G$ is denoted by $\delta(G)$ and $\Delta(G)$, respectively. A graph $G$ is $k$-\emph{regular} if every vertex has degree~$k$ in $G$. A $3$-regular graph is commonly referred to a \emph{cubic graph} in the literature. An \emph{isolated vertex} is a vertex of degree~$0$. If $X$ is a set of vertices of~$G$, then $\deg_X(v)$ is the number of neighbors of $v$ in $G$ that belong to the set $X$. In the special case when $X = V(G)$, we note that $\deg_X(v) = \deg_G(v)$. An $F$-\emph{component} of $G$ is a component of $G$ isomorphic to~$F$. A path and cycle on $n$ vertices is given by $P_n$ and $C_n$, respectively.

\section{Motivation and known results}
\label{S:known}

In 1998 Haynes and Slater~\cite{HaSl-98} proved that if $G$ is a connected graph of order $n \ge 3$, then $\gpr(G) \le n-1$ and they characterized the extremal family of graphs achieving equality in this bound. Moreover, Haynes and Slater~\cite{HaSl-98} proved that if $G$ is a connected graph of order $n$ and $\delta(G) \ge 2$, then $\gpr(G) \le \frac{2}{3}n$. Their proof contained an error, which was subsequently corrected by Huang and Shan~\cite{HuSh-11}. The graphs achieving equality in this $\frac{2}{3}n$-bound were characterized in~\cite{He-07}.
In 2007 Chen, Sun, and Xing~\cite{ChSuXi-07} proved that if $G$ is a cubic graph of order~$n$, then $\gpr(G) \le \frac{3}{5}n$. Goddard and Henning~\cite{GoHe-09} in 2009 showed that the only connected cubic graph achieving equality in this~$\frac{3}{5}n$-bound is the Petersen graph, and conjectured that if we exclude this exceptional graph, then the bound can be improved to $\gpr(G) \le \frac{4}{7}n$. Lu, Wang, Wang, and Wu~\cite{LuWaWaWu-12} in 2019 proved the conjecture in the special case of claw-free graphs. However, the conjecture has yet to be resolved in general.
In 2022 Henning, Pil\'{s}niak, and Tumidajewicz~\cite{HePiTu-22} proved that if $G$ is a graph of order~$n$ and $\delta(G) \ge 3$, then $\gpr(G) < \frac{127}{200} \, n = 0.635 \, n$. We summarize the best known upper bounds to date on the paired domination number for graphs with small minimum degree in Table~\ref{table1}.

\begin{center}
\begin{table}[htb]
{\small
\[
\begin{array}{||cccrcll||} \hline \hline %
& & & & & & \\
1998: & \delta(G) \ge 1 & \Rightarrow & \gamma(G) & \le &
 \displaystyle{ n-1  } & {\rm (\cite{HaSl-98})} \\
 & & & & & & \\
1998: & \delta(G) \ge 2 & \Rightarrow & \gamma(G) & \le & \displaystyle{  \frac{2}{3} \, n  }  & {\rm (\cite{HaSl-98,HuSh-11})} \\
& & & & & & \\
2022: & \delta(G) \ge 3 & \Rightarrow & \gamma(G) &  < & \displaystyle{ \frac{127}{200} \, n = 0.635 \, n \hspace*{0.15cm}^a } & {\rm (\cite{HePiTu-22})} \\
& & & & & & \\
\hline \hline
\end{array}
\]
\hspace*{2cm} $^a$ The precise bound given in~\cite{HePiTu-22} is $\gpr(G) \le \frac{19037}{30000}n < 0.634567 \, n$.
}
\begin{center}
\caption{Best known upper bounds on $\gpr(G)$ with small minimum degree $\delta \in [3]$.}
\label{table1}
\end{center}
\end{table}
\end{center}

For $k \ge 1$, let $\cG_k$ denote the class of all connected graphs with minimum degree at least~$k$ containing no isolated edge. We note that $\cG_1$ is the class of all connected graphs of order at least~$3$. The following problem is posed in~\cite{He-22}

\begin{problem}{\rm (\cite{He-22})}
\label{problem:Pdom}
Determine or estimate the best possible constants $c_{\pdom,k}$ (which depends only on~$k$) such that $\gpr(G) \le c_{\pdom,k} \cdot n(G)$ for all $G \in \cG_k$. These constants are given by
\[
c_{\pdom,k} =  \sup_{G \in \cG_k} \; \frac{\gpr(G)}{n(G)}.
\]
\end{problem}

The constants $c_{\pdom,k}$ are surprisingly only known for $k = 1$ and $k=2$. We summarize the best known results to date on the constants $c_{\pdom,k}$ for small $k \in [3]$ in Table~\ref{table2}.

\begin{center}
\begin{table}[htb]
{\small
\[
\begin{array}{||rcccll||} \hline \hline %
& & & & & \\
& & c_{\pdom,1} & = & 1  &  {\rm (\cite{HaSl-98})} \\
& & & & & \\
& & c_{\pdom,2} & = & \frac{2}{3} & {\rm (\cite{HaSl-98,He-07})} \\
& & & & & \\
\frac{3}{5} & \le & c_{\pdom,3} & < & \frac{3}{5} + \frac{13}{375} & {\rm (\cite{ChSuXi-07,GoHe-09,HePiTu-22})} \\
& & & & & \\
\hline \hline
\end{array}
\]
}
 \begin{center}
\caption{Best known upper bounds on $c_{\pdom,k}$ for small $k \in [3]$}
\label{table2}
\end{center}
\end{table}
\end{center}

\section{Main result}

It remains an open problem to determine the best possible upper bound on the paired domination number of a connected graph with minimum degree at least~$4$ in terms of its order~$n$. By the result in~\cite{HePiTu-22}, we have $\gpr(G) < 0.634567 \, n$. In this paper, we present an improvement of this best known upper bound to date on the paired domination number of a connected graph with minimum degree at least~$4$. We shall prove the following result, a proof of which is given in Section~\ref{S:proof}.

\begin{theorem}
\label{thm:main}
If $G$ is a graph of order~$n$ with $\delta(G) \ge 4$, then $\gpr(G) \le \frac{10}{17}n < 0.5883 n$.
\end{theorem}

\begin{figure}[htb]
\begin{center}
\begin{tikzpicture}[scale=0.85,style=thick,x=0.85cm,y=0.85cm]
\def\vr{2.5pt}
\path (0,0) coordinate (x1);
\path (3,0) coordinate (x2);
\path (3,3) coordinate (x3);
\path (0,3) coordinate (x4);
\path (1,1) coordinate (y1);
\path (1.1,1) coordinate (y1p);
\path (2,1) coordinate (y2);
\path (2,2) coordinate (y3);
\path (1,2) coordinate (y4);
\draw (x1) -- (x2) -- (x3) -- (x4) -- (x1);
\draw (y1) -- (y2) -- (y3) -- (y4) -- (y1);
\draw (x1) -- (y1) -- (x2);
\draw (x2) -- (y2) -- (x3);
\draw (x3) -- (y3) -- (x4);
\draw (x4) -- (y4) -- (x1);
\draw (x1) [fill=white] circle (\vr);
\draw (x2) [fill=black] circle (\vr);
\draw (x3) [fill=white] circle (\vr);
\draw (x4) [fill=black] circle (\vr);
\draw (y1) [fill=white] circle (\vr);
\draw (y2) [fill=black] circle (\vr);
\draw (y3) [fill=white] circle (\vr);
\draw (y4) [fill=black] circle (\vr);
\draw[anchor = north] (x1) node {{\small $x_1$}};
\draw[anchor = north] (x2) node {{\small $x_2$}};
\draw[anchor = south] (x3) node {{\small $x_3$}};
\draw[anchor = south] (x4) node {{\small $x_4$}};
\draw[anchor = north] (y1p) node {{\small $y_1$}};
\draw[anchor = west] (y2) node {{\small $y_2$}};
\draw[anchor = south] (y3) node {{\small $y_3$}};
\draw[anchor = east] (y4) node {{\small $y_4$}};
%
\end{tikzpicture}
\caption{The graph $H_8$}
\label{fig:H8H16}
\end{center}
\end{figure}
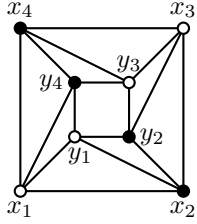

In the $4$-regular graph $G = H_8$ of order~$n = 8$ illustrated in Fig.\ref{fig:H8H16}, every edge is incident with a triangle. Therefore, $|N_G[u] \cup N_G[v]| \le 7$ holds for any two adjacent vertices $u$ and $v$. Thus, two adjacent vertices cannot form a dominating set of $G$, and so $\gpr(G) > 2$. Since the paired domination number of a graph is an even integer, this shows that $\gpr(G) \ge 4$. The set $\{x_2,x_4,y_2,y_4\}$ (indicated by the shaded vertices in Fig.\ref{fig:H8H16}) is an example of a PD-set of $G$ where $x_2$ and $y_2$ are paired and $x_4$ and $y_4$ are paired, and so $\gpr(G) \le 4$. Consequently, $\gpr(G) = 4 = \frac{1}{2}n$. As a consequence of Theorem~\ref{thm:main}, this yields the following bounds on $c_{\pdom,4}$.

\begin{corollary}
\label{cor:main}
$\frac{1}{2} \le c_{\pdom,4} \le \frac{10}{17} = \frac{1}{2} + \frac{3}{34}$.
\end{corollary}

\section{Proof of Theorem~\ref{thm:main}}
\label{S:proof}

We define the \emph{boundary} $\partial_G(S)$ of a set $S \subseteq V(G)$ in a graph $G$ as all neighbors of vertices of $S$ that belong outside the set $S$, that is, $\partial_G(S) = N_G[S] \setminus S$. We present the definition of a colored graph according to~\cite{HePiTu-22} and remark that this notion has a similar flavor to that of a residual graph defined in~\cite{Bu-21} (and also in other papers, such as in~\cite{Bu-22,BuHe-21,BuHe-25,BuKl-16}).

\begin{definition}{\rm (\cite{HePiTu-22})}
\label{residual-graph}
{\rm
Let $G$ be a graph and let $S$ be a set of vertices such that $G[S]$ contains a perfect matching. The \emph{colored graph} $G_S$ of $G$ associated with the set $S$ is the graph obtained from $G$ as follow: \\ [-21pt]
\begin{enumerate}
\item A vertex is colored \textbf{amber} if it has no neighbor in $S$.
\item A vertex is colored \textbf{beige} if it has a neighbor in $S$ and a neighbor not dominated by $S$.
\item A vertex is colored \textbf{cyan} if it and all its neighbors are dominated by $S$.
\item All edges not incident to amber vertices are removed from $G_S$. Thus every edge in $G_S$ joins two amber vertices or a beige and an amber vertex.
\end{enumerate}
}
\end{definition}

By Definition~\ref{residual-graph},  each vertex in the colored graph $G_S$ is colored amber, beige or cyan. In particular, a vertex in $S$ is colored cyan. We let $A$, $B$, and $C$ be the set of amber, beige, and cyan vertices, respectively, in $G_S$, and so $(A,B,C)$ is a partition of $V(G)$. The \emph{amber graph} is defined in~\cite{HePiTu-22} as the graph $G[A]$ induced by the set $A$ of amber vertices. The number of amber and beige vertices adjacent to a vertex $v$ in $G_S$ is the \emph{amber}-\emph{degree} and \emph{beige}-\emph{degree}, respectively, of~$v$, and is denoted by $\deg_A(v)$ and $\deg_B(v)$, respectively. The maximum amber-degree of a vertex in $A$ (resp., $B$) is denoted by $\Delta_A(A)$ (resp., $\Delta_A(B)$). If $v$ is an amber vertex, then its amber and beige neighbors are given by $N_A(v)$ and $N_B(v)$, respectively. We let $N_A[v] = N_A(v) \cup \{v\}$.

As observed in~\cite{HePiTu-22},  an amber vertex has no cyan neighbor, and therefore its degree in $G$ is the sum of its amber-degree and beige-degree in the colored graph $G_S$. Hence, the number of amber and beige neighbors of an amber vertex~$v$ in $G_S$ is precisely its degree in $G$, namely~$\deg_G(v) \ge \delta(G) \ge 4$. By construction of the colored graph, a beige vertex has at least one amber neighbor and possibly beige neighbors but has no cyan neighbor in $G_S$. Moreover, if $v$ is a beige vertex in $G_S$, then at least one of its neighbors in $G$ is colored cyan in $G_S$.

We are now in a position to prove Theorem~\ref{thm:main}. Recall its statement.

\newpage
\noindent \textbf{Theorem~\ref{thm:main}} \emph{If $G$ is a graph of order~$n$ with $\delta(G) \ge 4$, then $\gpr(G) \le \frac{10}{17}n < 0.5883 n$.
}

\noindent
\proof
If $\delta(G)\ge 4$ and $G$ contains an edge $e=uv$ such that $\deg_G(u) \ge \deg_G(v) \ge 5$, we delete this edge and obtain the graph $G'=G-e$. Then, the condition $\delta(G')\ge 4$ remains valid and $\gpr(G') \ge \gpr(G)$ also holds. Sequentially deleting edges between vertices of degree higher than~$4$, we obtain a graph $G''$ with $\delta(G'')=4$ and $\gamma(G'') \ge \gamma(G)$ such that  $\{v \in V(G'') \, \colon \deg_{G''}(v) \ge 5\}$ is an independent vertex set in $G''$. We may assume in the proof that $G$ already satisfies these properties. Since the paired domination number of a disconnected graph is the sum of the paired domination numbers of its components, we will also assume that $G$ is connected.
	
Let $G$ be a connected graph of order~$n$ such that $\delta(G) = 4$ and every edge is incident with at least one vertex of degree~$4$. Let $D$ be a subset of $V(G)$. We define $B_{\ge 4}$ as the set of beige vertices of degree at least~$4$ in $G_D$, and we define $B_i$ as the set of beige vertices of degree exactly~$i$ in $G_D$ for $i \in [4]$. Thus, every (beige) vertex in $B_i$ has exactly~$i$ amber neighbors for $i \in [4]$.  The weight $\w(v)$ of a vertex~$v$ is defined according to Table~3.

{\small
\[
\begin{array}{||c||c|c|c|c|c|c||} \hline
\mbox{set containing $v$} & A & B_{\ge 4} & B_3 & B_2 & B_1 & C \\ \hline
& & & & & &  \\
\w(v) & 45 & 33 & 31 & 29 & 27 & 0 \\
& & & & & & \\ \hline
\end{array}
\]
}
\begin{center}
\vskip -0.25 cm
\textbf{Table~3.} The weight $\w(v)$ of a vertex $v$
\end{center}

\noindent
We define the weight of the colored graph $G_D$ as the sum of the weights of its vertices, and so
\[
\w(G_D) = \sum_{v \in V(G)} \w(v) = 45|A| + 33|B_{\ge 4}| + 31|B_3| + 29|B_2| + 27|B_1|.
\]
\indent
The set $D$ is a PD-set in $G$ if and only if all vertices are colored cyan in $G_D$, in which case $\w(G_D) = 0$. Given the set $D \subseteq V(G)$ such that $G[D]$ contains a perfect matching and a set $S \subseteq V(G) \setminus D$ where $G[S]$ contains a perfect matching (and so, $|S| \ge 2$ is even), we define
\[
\xi(S) = \w(G_D) - \w(G_{D \cup S})
\]
as the decrease of the weight when extending $D$ to $D \cup S$. Such a set $S$ is a \emph{$D$-desirable set} if
\[
\xi(S) \ge 76.5|S|.
\]

If $D$ and $S$ are clear from the context, we set $D' = D \cup S$ and denote by $A'$, $B'$ and $C'$ the set of amber, beige, and cyan vertices, respectively, in $G_{D'}$. We define $B_i'$ as the set of beige vertices of degree exactly~$i$ in $G_{D'}$ for $i \in [4]$ and where $B_{\ge 4}'$ is the set of beige vertices of degree at least~$4$ in $G_{D'}$. The change of weight from an amber vertex $v \in A$ in $G_D$ to a beige vertex in $B_{\ge 4}'$ or to beige vertex in $B_i'$ in $G_{D'}$ where $i \in [3]$ is given in Table~4.

{\small
\[
\begin{array}{||c|c|c|c||} \hline
\mbox{$A  \rightarrow B_{\ge 4}$} & \mbox{$A \rightarrow B_3$} & \mbox{$A \rightarrow B_2$} & \mbox{$A \rightarrow B_1$} \\ \hline
& & & \\
12 & 14 & 16 & 18 \\
& & &\\ \hline
\end{array}
\]
}
\begin{center}
\vskip -0.25 cm
\textbf{Table~4.} Change of weight from an amber vertex $A$ to a beige vertex $B_i$
\end{center}

We note that the weight of an amber vertex in $G$ decreases by at least~$12$ when recolored beige in~$G_{D'}$. Next, we prove our key claim that if the weight of the colored graph $G_D$ is positive, then there exists a $D$-desirable set.

\begin{fact}
\label{f:fact1}
If $\w(G_D) > 0$, then the graph $G$ contains a $D$-desirable set.
\end{fact}
\proof Suppose, to the contrary, that $\w(G_D) > 0$, but the graph $G$ does not contain a $D$-desirable set with the given property. Suppose that $e$ is an edge incident with $y \in B_i$ in $G_D$ for some $i \in [4]$ and the edge $e$ joins $y$ to an amber vertex. If the edge $e$ is removed from $G_D$, then $y$ has exactly $i-1$ amber neighbors, implying that either $y$ is recolored cyan or $y$ belongs to the set $B_{i-1}$ in the colored graph $G_D - e$. By definition of the weights of the vertices given in Table~3, this in turn implies that the removal of the edge $e$ decreases the weight of $y$ by~$2$ if $i$ equals~$4$,~$3$, or~$2$ and by~$27$ if $i=1$. The change of weight from a beige vertex in $B_{i}$ to a beige vertex in $B_{i-1}$ where $i \in [4] \setminus \{1\}$, and from a beige vertex in $B_1$ to a cyan vertex is given in Table~5. We will frequently refer to these observations in the proof.

{\small
\[
\begin{array}{||c|c|c|c||} \hline
\mbox{$B_4 \rightarrow B_3$} & \mbox{$B_3 \rightarrow B_2$} & \mbox{$B_2 \rightarrow B_1$} & \mbox{$B_1 \rightarrow C$} \\ \hline
& & & \\
2 & 2 & 2 & 27 \\
& & & \\ \hline
\end{array}
\]
}
\begin{center}
\vskip -0.25 cm
\textbf{Table~5.} Change of weight from a beige vertex in $B_{i+1}$ to a beige vertex in $B_{i}$
\end{center}

In the proof, we present a series of subclaims stating structural properties of $G_D$ which culminate in the implication of its nonexistence. Throughout the proof of the claim, we let $v$ be an amber vertex of maximum amber-degree, and we let $w$ be a beige vertex of maximum amber-degree. Thus, $\deg_A(v) = \Delta_A(A)$ and $\deg_A(w) = \Delta_A(B)$.

Let $X = \partial(N_A[v])$ be the boundary of the set $N_A[v]$ in the amber graph $G[A]$, and so $X$ is the set of amber vertices that do not belong to $N_A[v]$ but have a neighbor in $N_A(v)$. Recall that if $x \in V(G)$, then $\deg_{X}(x)$ is the number of neighbors of $x$ that belong to the set $X$. Among all amber neighbors of~$v$, let $v'$ be chosen so that $\deg_{X}(v')$ is a maximum.

Let $Y = \partial(N_A[w])$ be the boundary of the set $N_A[w]$ in the amber graph $G[A]$, and so $Y$ is the set of amber vertices that do not belong to $N_A[w]$ but have a neighbor in $N_A(w)$. Among all amber neighbors of $w$, let $w'$ be chosen so that $\deg_{Y}(w')$ is a maximum.

\begin{claim}
\label{c:claim1}
$\Delta_A(A) \le 4$.
\end{claim}
\proof Suppose, to the contrary, that $\Delta_A(A) \ge 5$, and so $\deg_A(v) \ge 5$. By the edge-minimality of $G$, every neighbor of~$v$ has degree of at most~$4$ in $G$. Let $z$ be an arbitrary amber neighbor of $v$.

Suppose first that $\deg_{X}(z) = 0$, and so $N_A[z] \subset N_A[v]$. Thus, all amber neighbors of $z$ different from~$v$ are common amber neighbors of $z$ and $v$ in $G$. In this case, we let $z'$ be an amber neighbor of $v$ not in $N_A[z]$. We note that $z'$ exists since $\deg_A(v) \ge 5$ and $4 \ge \deg_A(z) $. We now let $S = \{v,z'\}$, and so $D' = D \cup S$. The amber vertices $v$, $z$ and $z'$ are colored cyan in the colored graph $G_{D'}$, while the further amber neighbors of $v$ are recolored as beige vertices and therefore belong to the set $B_i'$ in $G_{D'}$ for some $i \in [3]$ or are recolored as cyan vertices. We infer that each such neighbor of $v$ has a weight decrease of at least~$14$ (see Table~4). Therefore the total weight resulting from the set $S = \{v,z'\}$ is at least~$3 \times 45 + 3 \times 14 = 177 > 153 = 76.5|S|$. Thus, the set $S$ is a $D$-desirable set, a contradiction.

Hence, $\deg_{X}(z) \ge 1$, and so $N_A[z] \nsubseteq N_A[v]$ and $z$ has at least one neighbor that belongs to the set $X$. We now let $S = \{v,z\}$. The amber vertices $v$ and $z$ are colored cyan in the colored graph $G_{D'}$. As before, every amber neighbor of~$v$ different from~$z$ is recolored as a beige vertex $z' \in B_i'$ in $G_{D'}$ for some $i \in [3]$ or is recolored cyan, and therefore has a weight decrease of at least~$14$, while every neighbor of~$z$ that belongs to~$X$ is recolored in $G_{D'}$ as a beige vertex that belongs to the set $B_{\ge 4}'$ or $B_i'$ for some $i \in [3]$ or is recolored cyan, and therefore has a weight decrease of at least~$12$ (see Table~4). We therefore infer that the weight decrease in the colored graph $G_{D'}$ is at least~$2 \times 45 + 4 \times 14 + 1 \times 12  = 158 > 153 = 76.5|S|$. Thus, the set $S$ is a $D$-desirable set, a contradiction.~\smallqed

\begin{claim}
\label{c:claim2}
$\Delta_A(A) \le 3$.
\end{claim}
\proof Suppose, to the contrary, that $\Delta_A(A) = 4$. Let $v \in A$ be a vertex with $\deg_A(v)=4$ and $v'$ be the amber neighbor of $v$ as defined at the beginning of the proof. We set $S = \{v,v'\}$. Suppose that $\deg_{X}(v') = 0$, and so $N_A[v'] \subseteq N_A[v]$. Thus by our choice of~$v'$, every amber neighbor of $v$ is recolored cyan in $G_{D'}$, resulting in a total weight decrease of at least~$5 \times 45 = 225 > 153 = 76.5|S|$, implying that the set $S$ is a $D$-desirable set, a contradiction. Hence, $\deg_{X}(v') \ge 1$.

Consider first the case when $\deg_{X}(v') = 1$ and define $S = \{v,v'\}$ as before.  Every neighbor of $v'$ that belongs to~$X$ is recolored as a beige vertex $u \in B_i'$ in $G_{D'}$ for some $i \in [3]$ or is recolored cyan, and therefore has a weight decrease of at least~$14$. By the choice of $v'$, $\deg_{X}(z) \le 1$ holds for every $z \in N_A(v)$. Consequently, each neighbor of $v$ different from $v'$ belongs to $B_1'$ or is recolored cyan in $G_{D'}$ and has a weight decrease of at least~$18$. We therefore infer that the weight decrease in the colored graph $G_{D'}$ is at least~$2 \times 45 + 3 \times 18 + 1 \times 14 = 158 > 76.5|S|$. Thus, the set $S$ is a $D$-desirable set, a contradiction.

Now suppose that $\deg_{X}(v') \ge 2$. Let $S=\{v,v'\}$. Every neighbor of $v'$ in $N_X(v)$ and every amber neighbor of $v$ different from $v'$ is recolored as a beige vertex $v \in B_i'$ in $G_{D'}$ for some $i \in [3]$ or recolored cyan. Hence, the weight decrease is at least~$14$ for each of these five (or six) vertices. The vertices in $S$ are cyan vertices in $G_{D'}$. We therefore conclude that the weight decrease in $G_{D'}$ is at least~$2 \times 45 + 5 \times 14 = 160 > 76.5|S|$. This gives the desired contradiction and finishes the proof of the claim $\Delta_A(A) \le 3$.~\smallqed

\begin{claim}
\label{c:claim3}
$\Delta_A(B) \le 4$.
\end{claim}
\proof Suppose, to the contrary, that $\Delta_A(B) \ge 5$, and so $\deg_A(w) = \Delta_A(B) \ge 5$ where recall that $w$ is a beige vertex of maximum amber-degree. By Claim~\ref{c:claim2}, every neighbor of $w$ has amber-degree at most~$3$. Let $w'$ be as defined earlier so $\deg_Y(w')$ is maximum in $N_A(w)$. We set $S = \{w,w'\}$. The vertices $w$ and $w'$ are recolored cyan in $G_{D'}$, resulting in a weight decrease of~$33+45$.

If $\deg_{Y}(w') = 0$, the same is true for every amber neighbor of $w$. Thus in this case the vertex $w$ and all its amber neighbors in $G_D$ are colored cyan in $G_{D'}$, resulting in a weight decrease of at least~$33 + \deg_A(w) \times 45 \ge 33 + 5 \times 45  =258 > 76.5|S|$. Thus, the set $S$ is a $D$-desirable set, a contradiction. 

If $\deg_{Y}(w') = 1$, then each amber neighbor of $w$ different from $w'$ in $G_D$ is in $B_1'$ or is a cyan vertex in $G_{D'}$. Hence, for each such neighbor, the weight decrease is at least $18$. An amber neighbor of $w'$ in $Y$ belongs to  $B_1' \cup B_2'$ in $G_{D'}$, or it is a cyan vertex. In either case, the weight of the vertex decreases by at least $16$. We infer that the weight decrease in $G_{D'}$ is at least~$33+45+4 \times 18 +16=166 >76.5|S|$, a contradiction.

If $\deg_{Y}(w') \ge 2$, then each amber neighbor of $w$ and $w'$ belongs to a set $B_i'$ for some $i \in [3]$ or is a cyan vertex in $G_{D'}$ and consequently, its weight decreases by at least $14$. There are at least six such vertices different from  $w'$, and so the weight decrease in $G_{D'}$ is at least~$33+45+6 \times 14=162 >76.5|S|$, a contradiction.
~\smallqed

\medskip
By Claim~\ref{c:claim2}, we have $\Delta_A(A) \le 3$, and by Claim~\ref{c:claim3}, we have $\Delta_A(B) \le 4$. Thus, $B = B_1 \cup B_2 \cup B_3 \cup B_4$, where we recall that if $w \in B_i$ for $i \in [4]$, then $\deg_A(w) = i$.

\begin{claim}
\label{c:claim4}
$\Delta_A(A) \le 2$.
\end{claim}
\proof Suppose, to the contrary, that $\Delta_A(A) = 3$. Let $S = \{v,v'\}$. Suppose first that $\deg_{X}(v') = 0$. By our choice of~$v'$, every amber neighbor of $v$ is recolored cyan in $G_{D'}$, resulting in a total weight decrease of at least~$4 \times 45 = 180 > 76.5|S|$, implying that the set $S$ is a $D$-desirable set, a contradiction.

If $\deg_{X}(v') = 1$, then $\deg_{X}(z) \le 1$ holds for every amber neighbor $z$ of $v$ that is not in $N_A[v']$. Thus, $z \in B_1'$ or $z$ is a cyan vertex in $G_{D'}$, and its weight decreases by at least~$18$. The neighbor of $v'$ that belongs to~$X$ is recolored as a beige vertex $v \in B_i'$ in $G_{D'}$ for some $i \in [2]$ or is recolored cyan, and therefore has a weight decrease of at least~$16$.   By our assumption $\Delta_A(A) = 3$, every amber vertex has at least one beige neighbor in $G_D$. In particular, every vertex in $N_A[v] \cup N_A[v']$ has at least one beige neighbor. Thus, $\ell \ge 5$ edges join vertices in $N_A[v] \cup N_A[v']$ to beige vertices. As all vertices in $N_A[v] \cup N_A[v']$ are recolored beige or cyan in $G_{D'}$, these $\ell$ edges are not present in $G_{D'}$. The removal of them decreases the weight of the beige neighbors of vertices in $N_A[v] \cup N_A[v']$ by at least~$2 \times \ell$ (see Table~5). If $\ell \ge 6$, then
\[
\xi(S)\ge 2 \times 45 + 2 \times 18 + 1 \times 16 + 2 \times \ell \ge 154> 153= 76.5|S|,
\]
a contradiction. If $\ell=5$, then every vertex $u \in N_A[v] \cup N_A[v']$ has exactly one beige neighbor and consequently, $\deg_A(u)=3$ holds. For a vertex $z \in N_A(v)$ we also know that $\deg_{X}(z) \le 1$, that is, at most one amber neighbor of $z$ is outside $N_A[v]$. Therefore, every vertex $z \in N_A(v)$ has at least one neighbor from $N_A(v)$. Since $|N_A(v)|=3$, the subgraph induced by $N_A(v)$ in $G_D$ contains at least two edges. Consequently, there is a vertex $z' \in N_A(v)$ with $\deg_{X}(z) = 0$. Such a vertex $z'$ is a cyan vertex in $G_{D'}$, and we now can estimate the weight decrease and get a contradiction by
\[
\xi(S) \ge 3 \times 45 + 1 \times 18 + 1 \times 16 + 2 \times 5 =179> 153= 76.5|S|.
\]

If $\deg_{X}(v') \ge 2$, then $N_A[v] \cup N_A[v']$ contains at least six amber vertices. In the colored graph $G_{D'}$, vertices $v$ and $v'$ are recolored as cyan vertices; the remaining at least four vertices are either cyan or belong to $B_1' \cup B_2'$. For the latter case, the weight decreases by at least $16$ for each such vertex in~$B_1' \cup B_2'$. In $G_D$, there are $\ell \ge 6$ edges joining vertices in $N_A[v] \cup N_A[v']$ to beige vertices. The removal of these edges decreases the weight by at least~$2 \times \ell \ge 2 \times 6 = 12$ (see Table~5). We therefore infer $\xi(S) \ge 2 \times 45 + 4 \times 16 + 12=166 > 76.5|S|$. Thus, the set $S$ is a $D$-desirable set, a contradiction.~\smallqed

\medskip
By Claim~\ref{c:claim4}, we have $\Delta_A(A) \le 2$. Recall that by Claim~\ref{c:claim3}, we have $\Delta_A(B) \le 4$.

\begin{claim}
\label{c:claim5}
$\Delta_A(B) \le 3$.
\end{claim}
\proof Suppose, to the contrary, that $\Delta_A(B) = 4$, and so $\deg_A(w) = 4$. Since $\Delta_A(A) \le 2$, every amber neighbor of $w$ has amber-degree at most~$2$. Let $S = \{w,w'\}$. The vertices $w$ and $w'$ are recolored cyan in $G_{D'}$, resulting in a weight decrease of~$33 + 45$. If $\deg_{Y}(w') = 0$, then the vertex~$w$ and all its amber neighbors in $G_D$ are colored cyan in $G_{D'}$, resulting in a weight decrease of at least~$33 + 4 \times 45 =213 > 76.5|S|$. Thus, the set $S$ is a $D$-desirable set, a contradiction. Hence, $\deg_{Y}(w') \ge 1$.

Every amber neighbor of $w'$ which belongs to~$Y$ is recolored in $G_{D'}$ as a beige vertex that belongs to the set $B_1$ or is recolored cyan, and therefore has a weight decrease of at least~$18$ (see Table~4). Furthermore, every amber neighbor of~$w$ different from $w'$ is recolored in $G_{D'}$ as a beige vertex that belongs to the set $B_i$ for some $i \in [2]$ or is recolored cyan, and therefore has a weight decrease of at least~$16$ (see Table~4). We note that $\ell \ge 6$ edges join vertices in $N_A[w'] \cup N_A(w)$ to beige vertices different from $w$. As all vertices in $N_A[w'] \cup N_A(w)$ are recolored beige or cyan in $G_{D'}$, these $\ell$ edges are not present in $G_{D'}$. The removal of them decreases the weight of the beige neighbors of vertices in $N_A[w'] \cup N_A(w)$ by at least~$2 \times \ell \ge 2 \times 6 = 12$ (see Table~5). We therefore infer that the weight decrease in the colored graph $G_{D'}$ is at least~$33 + 45 + 3 \times 16 + 18 + 12 = 156 >76.5|S|$. Thus, the set $S$ is a $D$-desirable set, a contradiction.~\smallqed

\medskip
By Claim~\ref{c:claim5}, we have $\Delta_A(B) \le 3$. As  Claim~\ref{c:claim4} states $\Delta_A(A) \le 2$, and so the graph $G[A]$ induced by the amber vertices in $G_D$ consists of path and cycle components.

\begin{claim}
\label{c:claim6}
$G[A]$ contains no component that is a path $P_k$ of order $k \ge 3$.
\end{claim}
\proof Suppose, to the contrary, that $G[A]$ contains a component that is a path $P \colon  v_1 v_2 \dots v_k$ for some $k \ge 3$. If $k=3$, then let $u$ be a beige neighbor of $v_2$ and set $S=\{v_2, u\}$. In $G_{D'}$, the four vertices $v_1, v_2, v_3, u$  are recolored cyan. Without counting the further weight decreases, we get~$\xi(S) \ge 3 \times 45 +27=162 > 76.5|S|$, a contradiction.

If $P$ is a path component on $k \ge 4$ vertices, choose $S=\{v_2, v_3\}$. Vertices $v_1,v_2,v_3$ are colored cyan in $G_{D'}$, while $v_4$ is either cyan or belongs to $B_1'$. Without counting the decreases in the weights of beige vertices, we obtain~$\xi(S) \ge 3 \times 45 +18 =153 =76.5|S|$. Thus, the set $S$ is a $D$-desirable set, a contradiction.~\smallqed

\begin{claim}
\label{c:claim7}
$G[A]$ contains no cycle component.
\end{claim}
\proof Suppose, to the contrary, that $G[A]$ contains a component that is a cycle $C \colon v_1 v_2 \dots v_k v_1$ for some $k \ge 3$.

If $k \ge 8$, let $S=\{v_2,v_3, v_6,v_7\}$. In $G_{D'}$, each of the vertices $v_2, \dots v_7$ is cyan. Further, $v_1$ and $v_2$ are either cyan or belong to $B_1'$. Without counting the decrease in the weights of the adjacent beige vertices, we get~$\xi(S) \ge 6 \times 45 +2 \times 18= 306  = 76.5 |S|$, a contradiction.

If $k =7$, let $S=\{v_2,v_3, v_6,v_7\}$. In $G_{D'}$, all vertices of $C$ are recolored as cyan vertices. Without counting the further weight decreases, we obtain~$\xi \ge 7 \times 45= 315 > 306 = 76.5 |S|$ in $G_{D'}$. Thus, the set $S$ is a $D$-desirable set, a contradiction.

If $k =6$, we distinguish two cases. If there is a beige vertex $z$ that has a neighbor, say $v_1$, from $V(C)$ and a neighbor $z'$ outside $V(C)$, then we consider the set $S=\{z,v_1,v_4,v_5\}$. In this case the six vertices from $V(C)$ together with the vertex~$z$ are recolored cyan in $G_{D'}$. The (amber) neighbor $z'$ of $z$ belongs to $B_1' \cup B_2'$ or is colored cyan in $G_{D'}$. The weight of $z'$ therefore decreases by at least $16$. Without counting further decreases in weights, we get that~$\xi(S) \ge 6 \times 45+29+ 16=315  > 306 = 76.5 |S|$, once again yielding a contradiction.

In the other cases, all beige neighbors of the vertices on $C$ have neighbors only from $V(C)$. Let $S=\{v_1,v_2,v_4,v_5\}$ and assume that the neighborhood of $V(C)$ contains $b_3$, $b_2$, and $b_1$ vertices from $B_3$, $B_2$, and $B_1$, respectively. There are $\ell \ge 12$ edges between $V(C)$ and these beige vertices, and $\ell= 3b_3 +2b_2+b_1$. In $G_{D'}$, all these beige vertices are recolored as cyan, which means the following decrease in the weight of the beige vertices:
$$31  b_3 + 29 b_2 +27 b_1 > 10(3b_3+2b_2+b_1)=10 \ell \ge 120.
$$
Together with the decrease in the weight of vertices from the cycle, we obtain that the weight of $G_{D}$ decreases by at least $6 \times 45 + 120 =390 > 306 = 76.5 |S|$, which is a contradiction.

If $k=5$, similarly to the case of a $6$-cycle, we consider two possibilities. Suppose firstly that there is a beige vertex $z$ that has a neighbor, say $v_1$, from $V(C)$ and a neighbor $z'$ outside $V(C)$, consider the set $S=\{z,z',v_3,v_4\}$. In $G_{D'}$, the vertices in $V(C) \cup \{z,z'\}$ are all cyan. There are $\ell \ge 9$ edges in $G_D$ between $V(C) \cup \{z'\}$ and beige vertices different from $z$. The removal of these edges decreases the weights by at least $2 \times \ell \ge 18$. The total decrease in the weight of the colored graph is at least~$6 \times 45+ 29+ \18 = 317 > 306 = 76.5 |S|$, which is again a contradiction.

Suppose now that all beige neighbors of the vertices on $C$ have neighbors only from $V(C)$. Let $S=\{v_1,v_2,v_4,v_5\}$ and use the notation introduced for $k=6$. Now, we have $\ell \ge 10$ edges between $V(C)$ and the beige vertices and $\ell= 3b_3 +2b_2+b_1$. In $G_{D'}$, all these vertices are cyan, and the decrease in the weight of the beige vertices is
$$31  b_3 + 29 b_2 +27 b_1 > 10(3b_3+2b_2+b_1)=10 \ell \ge 100.
$$
As all vertices from $V(C)$ are recolored as cyan in $G_{D'}$, the decrease in the weight of the colored graph is at least~$5\times 45 +100=325 > 306 = 76.5 |S|$, a contradiction.

If $k=4$, we simply set $S=\{v_1,v_2\}$. Then, all vertices in $V(C)$ are cyan vertices in $G_{D'}$. It shows that the weight decrease in the colored graph $G_{D'}$ is larger than~$4 \times 45=180 >76.5|S| $. This gives the desired contradiction.

Finally, if $k=3$, we choose $S=\{v_1,z\}$, where $z$ is an arbitrary beige neighbor of $v_1$. In the colored graph $G_{D'}$, vertices $z$, $v_1$, $v_2$, and $v_3$ are colored cyan. This fact itself proves that the weight decrease in the colored graph $G_{D'}$ is larger than~$3\times 45 +27=162 >76.5|S| $. This contradiction finishes the proof of the claim.~\smallqed
\medskip

By Claims~\ref{c:claim6} and \ref{c:claim7}, the graph $G[A]$ consists of isolated vertices and $P_2$-components. In particular, $\Delta_A(A) \le 1$. We show next that beige vertex has amber-degree at most~$2$.

\newpage
\begin{claim}
\label{c:claim8}
$\Delta_A(B) \le 2$.
\end{claim}
\proof Suppose, to the contrary, that $\Delta_A(B) = 3$. Thus, $w$ denotes a beige vertex of amber-degree~$3$ and $w'$ is an amber neighbor of $w$ with a maximum $\deg_Y(w')$. Let $S=\{w,w'\}$. If $\deg_Y(w')=0$, then the four vertices in $N[w]$ are recolored cyan in $G_{D'}$ and~$\xi(S) \ge 31+3 \times 45=166 > 76.5|S|$, a contradiction.

If $\deg_Y(w')>0$, it equals $1$ and $w'$ belongs to a $P_2$-component $w'w''$ in $G[A]$. In $G_{D'}$,  vertices $w$, $w'$, and $w''$ are all cyan. The other two (amber) neighbors of $w$ are either cyan vertices in $G_{D'}$ or belong to $B_1'$. This shows that the total decrease in the weight of the colored graph is greater than~$31+2\times 45+2 \times 18= 157 > 76.5|S|$, a contradiction.~\smallqed

\medskip
By Claim~\ref{c:claim8}, we have $\Delta_A(B) \le 4$. Thus, $B = B_1 \cup B_2$.

\begin{claim}
\label{c:claim9}
No beige vertex has a neighbor from a $P_1$-component and from a $P_2$-component of $G[A]$.
\end{claim}
\proof Suppose, to the contrary, that $w \in B_2$ with two neighbors $w_1$ and $w_2$ such that $w_1$ belongs to a $P_1$-component and $w_2w_2'$ is a path in a $P_2$-component in $G[A]$. Let $S=\{w,w_2\}$. Since $w,w_1,w_2, w_2'$ are all colored as cyan vertices in $G_{D'}$, the decrease in the weight of the colored graph is greater than~$3 \times 45 + 29=164> 76.5|S|$, a contradiction.~\smallqed

\begin{claim}
\label{c:claim10}
$B_1=\emptyset$.
\end{claim}
\proof
Suppose, to the contrary, that $B_1 \ne \emptyset$. Suppose first that a vertex $z \in B_1$ has an amber neighbor $z'$ from a $P_1$-component of $G[A]$. If every beige neighbor of $z'$ belongs to $B_1$, we simply take $S=\{z,z'\}$ and observe that the entire $N[z']$ is recolored cyan in $G_{D'}$. Consequently, $\xi(S) \ge 45+4 \times 27=153 = 76.5|S|$, a contradiction. Hence at least one neighbor of $z'$ belongs to~$B_2$.

If $z_2\in B_2$ is a neighbor of $z'$, let $z_2'$ be the other neighbor of $z_2$. By Claim~\ref{c:claim9}, vertex $z_2$ belongs to a $P_1$-component in $G[A]$. Let $S=\{z',z_2\}$. In $G_{D'}$, the amber vertices $z'$, $z_2'$ and the beige vertices $z$, $z_2$ are recolored as cyan vertices. There are at least five edges between $\{z', z_2'\} $ and beige vertices different from $z$ and $z_2$. Therefore, $\xi(S) \ge 2 \times 45 + 27+29+ 5\times 2 = 156 > 76.5|S|$. This contradiction shows that no vertex from $B_1$ is adjacent to a vertex that belongs to a $P_1$-component of $G[A]$.

Now, suppose that a vertex $z \in B_1$ has an amber neighbor from a $P_2$-component $z_1z_2$  of $G[A]$. If there is a beige vertex $z'$ that has a neighbor $z_i$ from $\{z_1,z_2\}$ and an amber neighbor $u \notin \{z_1,z_2\}$, we set $S=\{z',z_i\}$. In $G_{D'}$, vertices  $z_1, z_2, z, z'$ are all cyan, and $u$ belongs to $B_1'$. Therefore, $\xi(S) \ge 2 \times 45 + 27 +29+ 18 = 164 > 76.5|S|$, a contradiction. If there is no beige vertex in $N_B[\{z_1,z_2\}]$ with a neighbor different from $z_1$ or $z_2$, then let $b_1= |B_1 \cap  N_B[\{z_1,z_2\}]|$ and  $b_2= |B_2 \cap  N_B[\{z_1,z_2\}]|$. Since there are at least six edges between $\{z_1,z_2\}$ and $B_1 \cup B_2$, we have $b_1 +2b_2 \ge 6$. Choosing $S=\{z_1,z_2\}$, all vertices in $\{z_1, z_2\} \cap  N_B[\{z_1,z_2\}]$ are recolored cyan in $G_{D'}$ and we have
\[
\xi(S) = 2\times 45 + 27 b_1 +29 b_2 \ge 90 + 14.5 \times (b_1+2b_2) \ge 90 + 14.5 \times 6 = 177 > 76.5 |S|.
\]
This contradiction finishes the proof of the claim.~\smallqed

\medskip
By Claims~\ref{c:claim8} and~\ref{c:claim10}, all beige vertices belong to $B_2$, that is, $B = B_2$. By Claim~\ref{c:claim9}, no beige vertex connects $P_1$- and $P_2$-components from $G[A]$. It follows then that $G[A \cup B]$ consists of two types of components. In a \emph{type-$1$ component}, every amber vertex $x$ has $\deg_A(x)=0$ and they are connected via $B_2$-vertices. In a \emph{type-$2$ component}, every amber vertex $x$ has $\deg_A(x)=1$ and the beige neighbors belong to $B_2$. In the continuation of the proof, we point out that neither type-$1$ nor type-$2$ components exist in $G_D$.

\begin{claim}
\label{c:claim11}
There is no type-$1$ component in $G_D$.
\end{claim}
\proof
Suppose, to the contrary, that there is a type-$1$ component $Q^1$ in $G_D$. By our earlier observations, every beige vertex in $Q^1$ has both its neighbors in $P_1$-component in $G[A]$. Since every vertex of $Q^1$ is incident to at least two edges, there is a cycle $C^1$ in $Q^1$. Observe that amber and beige vertices alternate along the cycle. Let $C^1 \colon v_1u_1v_2u_2 \dots v_ku_kv_1$ where $k \ge 2$ and where $v_i \in A$ with $\deg_A(v_i)=0$ and $u_i \in B_2$ for every $i \in [k]$.

Suppose firstly that $k$ is even. In this case, define
\[
S = \bigcup_{i=1}^{\frac{k}{2}} \{v_{2i-1},u_{2i-1}\},
\]
and so, $|S|=k$ and $G[S]$ contains a perfect matching. In $G_{D'}$, all vertices of $C^1$ are colored cyan decreasing the weight by $45k + 29k = 74k$. There are at least $2k$ edges between the amber vertices on the cycle $C^1$ and beige vertices that are not included in the cycle. Their removal further decreases the weight by at least $2\times 2k=4k$. Therefore, the weight of the colored graph decreases by at least $74k+4k =78k > 76.5k = 76.5|S|$. Thus, the set $S$ is a $D$-desirable set, a contradiction.

Suppose next that $k \ge 3$ is odd. We consider two cases. Suppose firstly that there is a beige vertex~$z$ having a neighbor from $C^1$, say it is $v_k$, and another (amber) neighbor $z' \notin V(Q^1)$. In this case, we let
\[
S = \{z,z'\} \cup \left( \bigcup_{i=1}^{\frac{k-1}{2}} \{v_{2i-1},u_{2i-1}\} \right)
\]
and so, $|S|=k+1$ and $G[S]$ contains a perfect matching. In $G_{D'}$ the vertices from $Z= V(Q^1) \cup \{z,z'\}$ are recolored cyan, resulting in a decrease in the weight by $45(k+1)+ 29(k+1)= 74(k+1)$. There are at least $2k+2$ edges between $Z \cap A$ and $B \setminus Z$, the removal of which decreases the weights of beige vertices by at least $2(2k+2)$. The total decrease in the weight of $G_D$ is then at least $74(k+1)+ 2(2k+2)= 78(k+1) > 76.5(k+1) = 76.5|S|$,  a contradiction.

Suppose next that there is no beige vertex that connects $V(C^1) \cap A$ to other amber vertices. In this case the component $Q^1$ consists of $v_1, \dots ,v_k$ and their beige neighbors. Since each $v_i$ has at least four beige neighbors, and every beige vertex has exactly two amber neighbors, we have $|V(Q^1) \cap B| \ge 2k$. We now let
\[
S = \bigcup_{i=1}^{\frac{k+1}{2}} \{v_{2i-1},u_{2i-1}\}
\]
and so, $|S|=k+1$ and $G[S]$ contains a perfect matching. In $G_{D'}$, all vertices of $Q^1$ are colored cyan, resulting in
\[
\xi(S) \ge 45k +29 |V(Q^1) \cap B| \ge 45k +58k = 103k > 76.5(k+1) = 76.5|S|,
\]
as $k \ge 3$. Hence, $S$ is a $D$-desirable set, and this contradiction proves the claim.~\smallqed

\begin{claim}
\label{c:claim12}
There is no type-$2$ component in $G_D$.
\end{claim}
\proof Suppose that $Q^2$ is a type-$2$ component which consists of $2k$ amber vertices and their beige neighbors. Since each amber vertex has at least three beige neighbors, and every beige vertex has exactly two amber neighbors, we have $|V(Q^2) \cap B| \ge \frac{3}{2} \times 2k=3k$. Now, let $S= V(Q^2) \cap A$. Since $Q^2$ is of type-$2$, we note that $|S| = 2k$ and $G[S]$ contains a perfect matching. In $G_{D'}$, all vertices of $Q^2$ are colored cyan and we obtain
\[
\xi(S) \ge 45\times 2k +29 |V(Q^2) \cap B| \ge 90k +29 \times 3k =177k > 76.5\times 2k = 76.5|S|.
\]
Consequently, $S$ is a $D$-desirable set and this contradiction proves the claim.~\smallqed

\medskip
Claims~\ref{c:claim11} and~\ref{c:claim12} complete the proof of Fact~\ref{f:fact1}.~\smallqed

\medskip
We now return to the proof of Theorem~\ref{thm:main}. By Fact~\ref{f:fact1}, if $\w(G_D) > 0$, then there is a $D$-desirable set in the graph $G$. Let $D_0 = \emptyset$ and let $G_0 = G_{D_0}$, and so $G_0$ is the graph $G$ with all vertices colored amber. We note that $V(G_0) = A$ and $\w(G_0) = 45n$. By Fact~\ref{f:fact1}, there exists a $D_0$-desirable set $S_1$, and so letting $D_1 = D_0 \cup S_1 = S_1$ and $G_1 = G_{D_1}$, we have $\w(G_0) - \w(G_1) \ge 76.5|S_1|$, that is,
\[
\w(G_1) \le \w(G_0) - 76.5|S_1|.
\]

If $\w(G_1) > 0$, then there is a $D_1$-desirable set $S_2$ by Claim~\ref{c:claim1}, and so letting $D_2 = S_1 \cup S_2$ and $G_2 = G_{D_2}$, we have $\w(G_1) - \w(G_2) \ge 76.5|S_2|$, that is,
\[
\w(G_2) \le \w(G_1) - 76.5|S_2|.
\]

If $\w(G_2) > 0$, then we continue the process, thereby obtaining a sequence of colored graphs $G_0, G_1, \ldots, G_k$ and a PD-set $D = S_1 \cup \cdots \cup S_k$ of $G$ such that
\[
\begin{array}{lcl}
0 = \w(G_k) & \le & \w(G_{k-1}) - 76.5|S_k| \1 \\
& \le & \displaystyle{ \w(G_0) - 76.5\sum_{i=1}^k |S_i| } \1 \\
& = & 45n - 76.5|D|.
\end{array}
\]
Consequently,
\[
\gpr(G) \le |D| \le \frac{45}{76.5}n = \frac{10}{17}n < 0.5883 \,n, 
\]
completing the proof of Theorem~\ref{thm:main}.~\QED

\section{Closing comment and conjectures}

In this paper we continue the study of upper bounds on the paired domination of a graph with given minimum degree in terms of its order. We establish a best known upper bound to date on the paired domination of a graph with minimum degree at least~$4$. However, several open problems and conjectures remain.

In 2007 Chen, Sun, and Xing~\cite{ChSuXi-07} conjectured that if $G$ is a connected graph of order $n \ge 11$ with $\delta(G) \ge 3$, then $\gpr(G) \le \frac{4}{7}n$. A slightly stronger conjecture was posed in 2009 by Goddard and Henning~\cite{GoHe-09}, namely that if $G$ is a connected graph of order $n$ with $\delta(G) \ge 3$, then $\gpr(G) \le \frac{4}{7}n$, unless $G$ is the Petersen graph, in which case $\gpr(G) = \frac{3}{5}n$. These two conjectures have yet to be settled. In this paper, we have shown (see Theorem~\ref{thm:main}) that if $G$ is a graph of order~$n$ and $\delta(G) \ge 4$, then $\gpr(G) \le \frac{10}{17}n < \frac{3}{5}n$. We do not know if this $\frac{10}{17}n$-upper bound is achievable. Motivated by the $\frac{4}{7}n$-conjectures stated above in~\cite{ChSuXi-07} and~\cite{GoHe-09}, we pose the conjecture that the $\frac{4}{7}n$-upper bound for the paired domination number holds if the minimum degree is at least~$4$. We state the conjecture formally as follows.

\begin{conjecture}
\label{conj:paired-dom}
If $G$ is a graph of order~$n$ with $\delta(G) \ge 4$, then $\gpr(G) \le \frac{4}{7}n = \left( \frac{10}{17} - \frac{2}{119} \right) n$.
\end{conjecture}

In 2014 Desormeaux and Henning in~\cite{DeHe-14} posed the conjecture that if $G$ is a bipartite cubic graph of order~$n$, then $\gpr(G) \le \frac{1}{2}n$. This conjecture has yet to be settled. Motivated by this $\frac{1}{2}n$-conjecture in~\cite{DeHe-14}, we close with the following conjecture.

\begin{conjecture}{\rm (\cite{DeHe-14})}
\label{conj:Pdom-3}
If $G$ is a bipartite $4$-regular graph of order $n$, then $\gpr(G) \le \frac{1}{2}n$.
\end{conjecture}

\section*{Acknowledgments}

Csilla Bujt\'{a}s acknowledges the support from the Slovenian Research Agency (ARIS) under the grants P1-0297 and N1-0355. Research of the second author was supported in part by the University of Johannesburg.

\end{document}